\documentclass[12pt]{article}
\usepackage{epsfig}
\usepackage{amssymb}
\usepackage{amsfonts}
\usepackage{amsmath,amsthm}
\usepackage{graphicx}
\usepackage{graphics}
\numberwithin{equation}{section}

\theoremstyle{definition}
\numberwithin{equation}{section}
\newcommand{\A}{\mathcal{A}}
\newcommand{\X}{\mathcal{X}}
\def\To{\longrightarrow}
\def\wt{\widetilde}
\def\La{\Lambda}
\def\la{\lambda}
\def\al{\alpha}
\title{APPROXIMATE WEAK AMENABILITY OF BANACH ALGEBRAS}   
\author{G. H. Esslamzadeh\ and B. Shojaee}
\date{}
\begin{document}


\maketitle

\pagestyle{plain}

\renewcommand{\theequation}{\arabic{section}.\arabic{equation}}

\baselineskip=18pt

\begin{abstract} In this paper we deal with four generalized notions of amenability which are called approximate, approximate weak, approximate cyclic and  approximate $n$-weak amenability. The first two were introduced and studied by Ghahramani and Loy in [9]. We introduce the third and fourth ones and we show by means of some examples, their distinction with their classic analogs.  

Our main result is that under some mild conditions on a given Banach algebra $\A$, if its second dual $\A^{**}$  is $(2n-1)$-weakly [respectively approximately/ approximately weakly/ approximately $n$-weakly] amenable, then so is $\A$. Also if $\A$ is approximately  $(n+2)$-weakly amenable, then it is approximately  $n$-weakly amenable. Moreover we show the relationship between approximate trace extension property and approximate weak [respectively cyclic] amenability. This answers question 9.1 of [9] for approximate weak and cyclic amenability.
\end{abstract}

\maketitle

\section{INTRODUCTION}

The concept of amenability for Banach algebras was introduced and studied for the first time by B. E. Johnson in [13]. Since then several variants of this concept have appeared in the literature each, as a kind of cohomological triviality. Bade, Curtis and Dales [1] introduced the notion of weak amenability for commutative Banach algebras and investigated  weak amenability of Beurling and Lipshitz algebras. Johnson [14] extended this concept to the non-commutative case and showed that group algebras of certain locally compact groups are weakly amenable. He proved the later result for arbitrary locally compact groups in [15]. Gronbaek in [12] and references therein, investigated properties of weakly amenable Banach algebras. In particular he showed that weakly amenable Banach algebras are essential. Moreover he identified the relationship between hereditary properties of weak and cyclic amenability and the trace extension property.  Dales, Ghahramani and Gronbaek introduced the notion of $n$-weak amenability as a generalization of weak amenability in [2],  where the authors discussed various cases in which $n$-weak amenability implies $m$-weak amenability and some cases which the converse implications do not hold. In particular they showed that $n$-weak amenability implies $(n+2)$-weak amenability. 

In all of the above mentioned concepts, all bounded derivations from a given Banach algebra $\A$ into certain Banach $\A$-bimodules  are required to be exactly inner. Gourdeau provided the following characterization of amenability: A Banach algebra $\A$ is amenable if and only if any bounded derivation from $\A$ into any Banach $\A$-bimodule is approximately inner, or equivalently weakly approximately inner [11, Proposition 2.1]. Moreover he showed that if $\A^{**}$ is amenable, then $\A$ is amenable. Motivated by Gourdeau's result, Ghahramani and Loy [9] introduced several approximate notions of amenability by requiring that all bounded derivations from a given Banach algebra $\A$ into certain Banach $\A$-bimodules to be approximately inner. However in contrast to Gourdeau's result, they removed the boundedness assumption on the net of implementing elements.  In the same paper and the subsequent one [10], the authors showed the distinction between each of these concepts and the corresponding classical notions and investigated properties of algebras in each of these new classes. 

At the begining, Ghahramani and Loy ask which of the standard results on amenability work for the approximate concepts [See 9, page 233], a question which identified the main direction of [9, 10] and the present paper. Motivated by this question, we study further, approximate amenability and approximate weak amenability. Moreover we introduce and study two new concepts of, approximate $n$-weak amenability and approximate cyclic amenability. 
In question 9.1 of [9] the authors ask, What are the hereditary properties of approximate concepts of amenability? We answer this question for approximate weak and cyclic amenability by defining approximate trace extension property  and constracting approximate analogs of certain results of Gronbaek in [12]. Then we provide a characterization of approximate amenability. That a Banach algebra $\A$ inherits approximate amenability and approximate weak amenability from $\A^{**}$, is the main result of section 2. Indeed the first one appeared in [9], but here we provide a short proof for it. In Section 3, first we show that approximate $(n+2)$-weak amenability implies approximate $n$-weak amenability; Then we discuss the relationship between approximate $n$-weak amenability of $\A$ and its unitization. Passage of $n$-weak and approximate $n$-weak amenability from $\A^{**}$ to $\A$ for odd $n$, is the main result of section 3. In the last section we discuss hereditary properties of approximate cyclic amenability and the relationship between approximate cyclic amenability of $\A$ and its unitization. Some of our arguments were inspired by their classic analogs mostly in [2, 3, 9, 12].
 
Before proceeding further we recall some terminology.

Throughout $\A$ is a Banach algebra and $\X$ is a Banach $\A$-bimodule. Also $\A$-module means Banach $\A$-bimodule. A linear mapping $D:\A\To\X$ is a derivation if $D(ab)=a.Db+Da.b$ for $a,b\in \A$. For any $x\in \X$, the mapping $\delta_{x}:a\longmapsto ax-xa$ ,$a\in \A$, is a continuous derivation which is called an inner derivation. Let $Z^{1}(\A,\X)$ be the space of all continuous derivations from $\A$ into $\X$ and $B^{1}(\A,\X)$ be the space of all inner derivations from $\A$ into $\X$. Then the first cohomology group of $\A$ with coefficients in $\X$ is the quotient space $H^{1}(\A,\X)=Z^{1}(\A,\X)/B^{1}(\A,\X).$ For each $n\geq1$, $\A^{(n)}$, the n-th conjugate space of $\A$, is an $\A$-module, with the module actions defined inductively by
$$
\langle u,\La.a\rangle=\langle a.u,\La\rangle\quad ,\quad \langle u,a.\La\rangle=\langle u.a,\La\rangle\quad, \La\in \A^{(n)}, u\in \A^{(n-1)}, a\in \A. 
$$
 The {\it first} and {\it second Arens multiplications on $\A^{**}$} that we denote by ``$\Box$'' and ``.'' respectively, are defined in three steps. For $a,b\in\A,\ f\in\A^*$ and $m,n\in\A^{**}$, the elements $f\Box a,\ a.f,\ m\Box f,\ f.m$ of $\A^*$ and $m\Box n,\ m.n$ of $\A^{**}$ are defined in the following way:
$$
\alignedat2 &\langle f\Box a\ ,\ b\rangle =\langle f\ ,\ ab\rangle \qquad\qquad &&\langle a.f\ ,\ b\rangle =\langle f\ ,\ ba\rangle \\&      \langle m\Box f\ ,\ a\rangle =\langle m\ ,\ f\Box a\rangle &&\langle f.m\ ,\ b\rangle =\langle m\ ,\ b.f\rangle \\&   \langle m\Box n\ ,\ f\rangle =\langle m\ ,\ n\Box f\rangle &&\langle m.n\ ,\ f\rangle =\langle n\ ,\ f.m\rangle  .\endalignedat
$$ 
The second dual of a Banach algebra, equipped with the first [respectively second] Arens product is a Banach algebra. 
 We always consider the second dual of a Banach algebra with the first Arens product. 

A Banach algebra $\A$ is called $n$-weakly amenable if $H^{1}(\A,\A^{(n)})={0}.$ Usually, 1-weakly amenable Banach algebras are called weakly amenable. For a detailed discussion of 2-weak amenability of weighted group algebras the reader can see [5, Chapter 12]. A derivation $D:\A\To \X$ is called approximately inner if there exists a net $(x_{\al})\subseteq \X$ such that $D(a)=\lim_{\al}(a.x_{\al}-x_{\al}.a)\quad (a\in \A).$ We say that  $\A$ is  approximately amenable if for each Banach $\A$-module $\X$ every bounded derivation $D:\A\To \X^{*}$ is approximately inner and it is approximately weakly amenable if every bounded derivation $D:\A\To \A^{*}$ is approximately inner. 
Obviously  if  $\A$ is approximately amenable, then it is approximately weakly amenable. However the converse is not true. For instance,  let $F_{2}$ be the free group with two generators. It is known that $\ell^{1}(F_{2})$ is approximately weakly amenable but by [9, Theorem 3.2], it is not approximately amenable. If $\A$ is weakly amenable then it is approximately weakly amenable but the converse is not true as it was shown in [9, example 6.2]. In the case that $\A$ is commutative, the only inner derivation from $\A$ into $\A^*$ is the zero derivation and hence $\A$ is weakly amenable if and only if $\A$ is approximately weakly amenable.

\section{APPROXIMATE WEAK AMENABILITY}

Recall that $\A$ is called essential if $\overline{\A^{2}}=\A$. The following proposition can be proved  by using the arguments of [3, proposition 1.3] and [2, Theorem 2.8.63] with proper modifications.    

\textbf{Proposition 2.1.} Suppose $\A$ is approximately weakly amenable. Then,

(i) $\A$ is essential;

 If $\A$ is commutative, then 

(ii) $Z^{1}(\A,\X)=\{0\}$ for every $\A$-module $\X$.

(iii) if $\theta:\A\To \mathcal B$ is a continuous homomorphism into the Banach algebra $\mathcal B$ such that $\overline{\theta(\A)}=\mathcal B$,  then $\mathcal B$ is approximately weakly amenable.

\textbf{Definition.} Let $I$ be a closed ideal in $\A$. We say that $I$ has the approximate trace extension property if for each $\la\in I^{*}$ with  $a.\la=\la.a (a\in\A)$ there is a net $(\La_{\al})\subseteq\A^{*}$ such that
$$
\La_{\al}|_{I}=\la\ (\text{for\ any}\ \al)\quad\quad \text{and} \quad\quad a.\La_{\al}-\La_{\al}.a\To 0\quad (a\in\A).
$$
\textbf{Definition.} Let $I$ be a closed ideal in $\A$. We say that a bounded approximate identity 
$\{ e_\al\}$ of $I$ is  quasi central for $\A$ if $\lim_\al\|ae_\al -e_\al a\|=0$ for all $a\in\A$. 

\textbf{Proposition 2.2.} Let $I$ be a closed ideal in $\A$.

(i) Suppose that $\A/I$ is approximately weakly amenable. Then $I$ has the approximate trace extension property. 

(ii) Suppose that $\A$ is weakly amenable and $I$ has the approximate trace extension property.
Then $\A/I$ is approximately weakly amenable.

(iii) Suppose that $I$ is weakly amenable and $\A/I$ is approximately weakly amenable. Then $\A$ is approximately weakly amenable.

(iv) Suppose that $\A$ is approximately weakly amenable and $I$ has a quasi-central bounded approximate identity for $\A$. Then $I$ is approximately weakly amenable.

{\it Proof.} (i) Let $\la\in I^*$ is such that $a.\la=\la.a (a\in\A)$. Take $\La\in\A^{*}$ with $\La|_{I}=\la$. Define
$$D:a+I\longmapsto a.\La-\La.a,\quad\A/I\To I^{\perp}=(\A/I)^{*}.$$

We see immediately that $D\in Z^{1}(\A/I,(\A/I)^{*}).$ Since $\A/I$ is
approximately weakly amenable, there exists a net
$(\La_{\al})\in I^{\perp}$ such that
$$D(a+I)=\lim_{\al}(a.\La_{\al}-\La_{\al}.a)\quad (a\in\A).$$

Set $\tau_{\al}=\La-\La_{\al}\in\A^{*}.$ Then
$$\tau_{\al}|_{I}=\la\ (\text{for\ any}\ \al)\quad \text{and} \quad \lim_{\al}(a.\tau_{\al}-\tau_{\al}.a)=0\quad (a\in\A).$$

(ii) Suppose $\pi:A\To\A/I$ is the quotient map and
$D\in Z^{1}\in (\A/I,(\A/I)^{*})$. Set $\wt{D}=\pi^{*}\circ
D\circ\pi$. Then $\wt{D}\in Z^{1}(A,A^{*})$, and so there
exists $\La\in\A^{*}$ with
$$\wt{D}a=a.\La-\La.a \quad (a\in\A).$$

Clearly $\wt{D}a|_I=0$ and hence $a.\La|_I=\La|_I.a\ (a\in\A)$. Thus by assumption there exists a net $(\La_{\al})\subseteq\A^{*}$ such that
for every $\al$ we have $\La_{\al}|_{I}=\La|_{I}$ and
$$\lim_{\al}(a.\La_{\al}-\La_{\al}.a)=0\quad (a\in\A).$$

Then $\La-\La_{\al}\in I^{\perp}$ and
$$D(a+I)=\lim_{\al}[a.(\La-\La_{\al})-(\La-\La_{\al}).a]\quad (a\in\A).$$

Now it follows that $\A/I$ is approximately weakly amenable.

(iii) Suppose $\iota:I\To\A$ is the natural
embedding and $D\in Z^{1}(\A,\A^{*})$. Then $\iota^{*}\circ
D\circ\iota\in Z^{1}(I,I^{*})$, and so, since I is weakly amenable,
there exists $\La_{1}\in I^{*}$ with
$$
(\iota^{*}\circ D)(a)=\delta_{\La_1}(a)\quad (a\in I).
$$
Extend $\La_{1}$ to be an element of $A^{*}$ and call it again $\La_{1}$. By replacing $D$ by
$D-\delta_{\La_1}$, we suppose that $(\iota^{*}\circ D)|_{I}=0.$

For $a,b\in I$ and $c\in\A$, we have
$$\langle c,Dab\rangle=\langle ca,\iota\circ D(b)\rangle+\langle bc,\iota^{*}(a)\rangle=0,$$

and hence $D|_{I^{2}}=0$. By [2, Theorem 2.8.63(i)], $\overline{I^{2}}=I$, and so
$D|_{I}=0$. If we set $F=\overline{I\A +\A I}$, then $F=\overline{I^{2}}=I$. For each
$a\in\A$ and $b\in I$, we have $a.Db=Dab=0$, and so $Da.b=0$. Taking
$c\in\A$, we get
$$
\langle b.c,Da\rangle=\langle c,Da.b\rangle=0
$$
and so $Da|_{I\A}=0$. Similarly $Da|_{\A I}=0$, and hence $Da|_{I}=0.$
Thus $D(\A)\subseteq I^{\perp}$ and the map
$$
\wt{D}:\A/I\To I^{\perp},\quad a+I\longmapsto Da
$$
is a continuous derivation. By hypothesis $\A/I$ is approximately
weakly amenable, and so there exists a net
$(\La_{\al})\subseteq I^{\perp}$ such that
$$
Da=\lim_{\al}(a.\La_{\al}-\La_{\al}.a)\quad (a\in\A).
$$
Therefore
$$Da=\lim_{\al}[a.(\La_{\al}+\La_{1})-(\La_{\al}+\La_{1}).a]\quad(a\in\A).$$

It now follows that $\A$ is approximately weakly amenable. 

(iv) By [12, Proposition 1.3] any bounded derivation $D:I\To I^*$ can be lifted  to a bounded derivation $\wt D:\A\To \A^*$, from which, the result follows immediately.          $\qed$

In the following Proposition we characterize approximately  amenable Banach algebras. 

\textbf{Proposition 2.3.} The following conditions are equivalent;

\quad(i) $\A$ is approximately amenable;

\quad(ii) For any $\A$-module $\X$, every bounded derivation
$D:A\To \X^{**}$ is approximately inner.

{\it Proof.} (i)$\Longrightarrow$(ii) This is immediate.

(ii)$\Longrightarrow$(i) By [9, Theorem 2.1], it suffices to show that
if $D\in Z^{1}(\A,\X)$, then it is approximately inner.

We have $\iota\circ D\in Z^{1}(\A,\X^{**})$, where
$\iota:\X\To \X^{**}$ is the canonical embedding. By (ii),
there exists a net $(\La_{\al})\subseteq \X^{**}$ such that
$$\iota\circ D(a)=\lim_{\al}( a.\La_{\al}-\La_{\al}.a)\quad (a\in\A).$$

Now take $\epsilon>0$, and non-empty finite sets $E\subseteq\A$, $F\subseteq \X^{*}$. Then there is an $\al$ such that
$$
|\langle \la,\iota\circ D(a)-(a.\La_{\al}-\La_{\al}.a)\rangle|<\epsilon
$$
for all $\la\in F$ and $a\in E.$ By Goldstine's theorem there is an $x_{\al}\in \X$ such that
$$
|\langle \la,\iota\circ D(a)-(a.x_{\al}-x_{\al}.a)\rangle|<\epsilon
$$
for all $\la\in F$ and $a\in E.$ Thus there is a net $(x_{\al})\subseteq \X$ such that
$$
Da=\omega-\lim_{\al}(a.x_{\al}-x_{\al}.a)\quad (a\in\A).
$$
Finally, for each finite set $E\subseteq\A$, say
$E=\{a_{1},\ldots,a_{n}\}$,
$$
(a_{1}.x_{\al}-x_{\al}.a_{1},\ldots, a_{n}.x_{\al}-x_{\al}.a_{n})\To(Da_{1},\ldots,Da_{n})
$$
weakly in $\X^{n}$. By Mazur's Theorem,
$$
(Da_{1},\ldots,Da_{n})\in \overline{Co^{\|.\|}}\{(a_{1}.x_{\al}-x_{\al}.a_{1},\ldots, a_{n}.x_{\al}-x_{\al}.a_{n})\}.
$$
Thus there is a convex linear combination $x_{E,\epsilon}$ of
elements of the set $\{x_{\al}\}$ such that
$$
\|Da-(a.x_{(E,\epsilon)}-x_{(E,\epsilon)}.a)\|<\epsilon\quad (a\in E).
$$
The family of such pairs $(E,\epsilon)$ is a directed set for the
partial order $\leq$ given by 
$$
(E_{1},\epsilon_{1})\leq(E_{2},\epsilon_{2})\ \text{if} \ E_{1}\subseteq E_{2}\ \text{and}\ 
\epsilon_{1}\geq\epsilon_{2}
$$
and
$$
Da=\lim_{(E,\epsilon)}a.x_{(E,\epsilon)}-x_{(E,\epsilon)}.a\quad (a\in\A). \qed
$$ 

In the first part of the following Theorem, we provide a short proof for [9, Theorem 2.3]. The reader should be aware that we can turn $\A^{***}$ into an $\A^{**}$-module in two differnt ways. First, consider $\A^{**}$ as an $\A^{**}$-module and then equip $\A^{***}$ with the usual dual module structure. In the second method, $\A^*$ is considered as an $\A$-module and taking the second duals, $\A^{***}$ turns into an $\A^{**}$-module as in [11]. A similar statement holds for higher duals $\A^{(2n-1)}$ for $n\ge 2$. In the next Theorem and Theorem 3.6 we deal with the first one of these module structures. Also the term  $B(\A,\A^*)$ [respectively $W(\A,\A^*)$] denotes the space of bounded [respectively weakly compact] operators from $\A$ into $\A^*$.  

\textbf{Theorem 2.4.}(i) Suppose that $(\A^{**},\Box)$ is approximately amenable. Then $\A$
is approximately amenable.

(ii) Suppose that $B(\A,\A^*)=W(\A,\A^*)$ and $(\A^{**},\Box)$
is approximately weakly amenable. Then $\A$ is approximately weakly amenable.

{\it Proof.} (i) Let $D\in Z^{1}(\A,\X^{*})$ for an $\A$-module $\X$. As in
[11], $\X^{***}$ turns into an $(\A^{**},\Box)$-module so that
$D^{**}:(\A^{**},\Box)\To \X^{***}$
is a continuous derivation. By [10, Theorem 2.1] $\A^{**}$ is approximately contractible. So there exists a net
$(\La_{\al})\subseteq \X^{***}$ such that
$$
D^{**}\phi=\lim_{\al}(\phi.\La_{\al}-\La_{\al}.\phi) \quad (\phi\in\A^{**}).
$$
Let $P:\X^{***}\To \X^{*}$ be the natural projection. Then
$$
Da=P(Da)=\lim_{\al}(a.P(\La_{\al})-P(\La_{\al}).a)\quad (a\in\A),
$$
and so D is approximately inner. Therefore $\A$ is approximately amenable.

(ii) Let $D\in Z^{1}(\A,\A^{*})$. Consider $\A^{**}$ as an $\A^{**}$-module and let $\A^{***}$ be its dual module. As in the proof of [2, Proposition 2.8.59(iii)] using the assumption $B(\A,\A^*)=W(\A,\A^*)$ we see that  $D^{**}:(\A^{**},\Box)\To\A^{***}$
is a continuous derivation. Thus there exists a net $(\La_{\al})$ in  $\A^{***}$ such that
$$
D(\phi)=\lim_{\al}(\phi .\La_{\al}-\La_{\al}.\phi)\quad(\phi\in\A^{**}).
$$
Let $P:\A^{***}\To\A^{*}$ be the natural projection. Then
$$
Da=P(Da)=\lim_{\al}(a.P(\La_{\al})-P(\La_{\al}).a)\quad (a\in\A).
$$
Therefore $\A$ is approximately weakly amenable.$\qed$

\section{APPROXIMATE $N$-WEAK AMENABILITY}

\textbf{Definition.} Let $n$ be a natural number. We say that $\A$ is approximately $n$-weakly amenable if every continuous derivation $D:\A\To\A^{(n)}$ is approximately inner; $\A$ is permanently approximately weakly amenable if $\A$ is approximately
$n$-weakly amenable for each $n\in\Bbb N.$

The following observations are immediate consequences of the above definition.

(i) An approximately amenable Banach algebra is permanently approximately weakly amenable.

(ii) A commutative Banach algebra is $n$-weakly amenable if and only if it is approximately $n$-weakly amenable.

 (iii) A commutative Banach algebra is permanently approximately weakly
amenable if and only if it is approximately weakly amenable.

\textbf{Theorem 3.1.}  Suppose that $\A$ is approximately $(n+2)$-weakly amenable.
Then $\A$ is approximately $n$-weakly amenable.

{\it Proof.} Let $D\in Z^{1}(\A,\A^{n})$. Then D can be viewed as an element
of $Z^{1}(\A,\A^{(n+2)})$, and so there exists a net
$(\La_{\al})\in\A^{(n+2)}$ with
$$
Da=\lim_{\al}(a.\La_{\al}-\La_{\al}.a)\quad (a\in\A).
$$
Let $P:\A^{(n+2)}\To\A^{(n)}$ be the natural projection. Then
$$
Da=P(Da)=\lim_{\al}(a.P(\La_{\al})-P(\La_{\al}).a)\quad (a\in\A).
$$
So D is an approximately inner derivation. Thus $\A$ is 
approximately $n$-weakly amenable. $\qed$

\textbf{Remark 3.2.} Suppose $\A$ is not unital and $\A^{\#}=\Bbb C e\oplus\A$ is its unitization. Define $e^{*}\in\A^{\#^{*}}$ by requiring that $\langle
e^{*},e\rangle=1$ and $e^{*}|_{\A}=0$. Then for every $n\in\Bbb N$ we have the
identifications
$$\A^{\#^{(2n)}}=\Bbb C e\oplus\A^{(2n)},$$
$$\A^{\#^{(2n-1)}}=\Bbb C e^{*}\oplus\A^{(2n-1)}.$$
The module actions of $\A^{\#}$ on $\A^{\#^{(2n-1)}}$ are given by
$$
(re+a).(se^{*}+\La)=(rs+\langle \La,a\rangle)e^{*}+r\La+a.\La ,
$$
$$
(se^{*}+\La).(re+a)=(rs+\langle \La,a\rangle)e^{*}+r\La+\La.a.
$$
where $r,s\in \Bbb C $, $a\in\A^{\#}$ and $\La\in\A^{\#^{(2n-1)}}$. The assymetry of results for odd and even duals in the next proposition is due to the fact  that $\A^{(2n-1)}$ is not a submodule of $\A^{\#^{(2n-1)}}$ in general, but $\A^{(2n)}$ is a submodule of $\A^{\#^{(2n)}}$.

\textbf{Proposition 3.3.} Let $\A$ be a non-unital Banach algebra and $n\in\Bbb N$. 

(i) If $\A^{\#}$ is approximately $2n$-weakly amenable, then
$\A$ is approximately $2n$-weakly amenable.

(ii) If $\A$ is approximately $(2n-1)$-weakly amenable, then
$\A^{\#} $ is approximately $(2n-1)$-weakly amenable.

(iii) Suppose that $\A$ is commutative. Then $\A^{\#}$ is
approximately $n$-weakly amenable if and only if $\A$ is approximately
$n$-weakly amenable.

(iv) Suppose that $\A$ has a bounded approximate identity. Then $\A$ is approximately weakly amenable if and only if $\A^{\#}$ is approximately weakly amenable.

{\it Proof.} (i) Let $D\in Z^{1}(\A,\A^{2n})$ and define
$$
\wt{D}:\A^{\#}\To\A^{2n},\quad re+a\longmapsto Da,\quad r\in\Bbb C,\quad a\in\A .
$$ 
Certainly $\wt{D}$ is a continuous derivation and hence can be viewed as an
element of $Z^{1}(\A^{\#},\A^{\#(2n)})$. So there exists a net
$(\La_{\al}^{'})\subseteq\A^{\#(2n)}$ such that
$$
Da=\lim_{\al}(a.\La_{\al}^{'}-\La_{\al}^{'}.a)\quad (a\in\A).
$$
Using the identity $\A^{\#(2n)}=\Bbb C e^*\oplus\A^{(2n)}$ for every $\al$ we find  $\La_{\al}\in\A^{(2n)}$ and $r_{\al}\in \Bbb C $ such that $\La_{\al}^{'}=r_{\al} e^*+\La_{\al}$  and so
$$
Da=\lim_{\al}(a.\La_{\al}-\La_{\al}.a).
$$
(ii) Let $D:\A^{\#}\To\A^{\#(2n-1)}$ be a continuous derivation. Then $D|_{\A}$ which we denote by $D$ again, is a continuous derivation and by Remark 3.2, is of the form
$$
D:\A\To \Bbb C  e^{*}\oplus\A^{2n-1},\quad a\longmapsto \langle a,\la\rangle e^{*}+\wt{D}(a).
$$
 It is easy to see that $\wt{D}:\A\To\A^{2n-1}$ is a continuous derivation,
and so there exists a net $(\La_{\al})\subseteq\A^{2n-1}$
such that
$$
\wt{D}(a)=\lim_{\al}(a.\La_{\al}-\La_{\al}.a)\quad (a\in\A).
$$
Let $a,b\in\A$. Then we have
$$
\aligned\langle \la,ab\rangle&=\langle\wt{D}b,a\rangle+\langle \wt{D}a,b\rangle\\
&=\lim_{\al}\langle b.\La_{\al}-\La_{\al}.b,a\rangle
+\lim_{\al}\langle a.\La_{\al}-\La_{\al}.a,b\rangle\\
&=\lim_{\al}(\langle \La_{\al},ab\rangle-\langle \La_{\al},ba\rangle
+\langle \La_{\al},ba\rangle-\langle \La_{\al},ab\rangle)=0. \endaligned
$$
So $\la|_{\A^{2}}=0$. On the other hand by Theorem 3.1 $\A$ is
approximately weakly amenable, and hence by Proposition 2.1(i),
$\A^{2}$ is dense in$\A$. It follows that $\la=0$, and so
$D=\wt{D}$ is an approximately inner derivation.

(iii) We observed at the begining of this section that for commutative Banach algebras approximate $n$-weak amenability and $n$-weak amenability are equivalent. So this part follows from [3, Proposition 1.4 (iii)]. 

(iv) Part (ii) implies that, if $\A$ is approximately weakly amenable, then so is $\A^{\#}$. Conversely suppose  $\A^{\#}$ is approximately weakly amenable and $D:\A\To\A^*$ is a continuous derivation. By [12, Corollary 2.2], $D$ can be extended to a continuous derivation $\wt{D}:\A^{\#}\To(\A^{\#})^*$. By assumption    and Remark 3.2 there is a net $\{\La_{\al}\}$ in $\A^*$ and a net $\{ s_{\al}\}$ in $\Bbb C$ such that for every $a\in\A$ and $r\in\Bbb C$ 
$$
\wt{D}(re+a)=\lim_{\al}[(re+a)(s_{\al}e^*+\La_{\al})-(s_{\al}e^*+\La_{\al})(re+a)]=\lim_{\al}[a.\La_{\al}-\La_{\al}.a]. 
$$
Thus 
$$
D(a)=\lim_{\al}[a.\La_{\al}-\La_{\al}.a]\quad a\in\A .
$$
Therefore $\A$ is approximately weakly amenable. $\qed$

\textbf{Remark 3.4.} From part (ii) of the preceding proposition we conclude that if $\A$ approximately weakly amenable then so is $\A^{\#}$; But we do not know whether its converse is true in general. However in the case that $\A$ has a bounded approximate identity the converse is true as we saw in part (iv) of the preceding proposition. For approximate cyclic amenability  the situation is different as we see in the Proposition 4.1.

 The following theorem is a partial converse of Theorem 3.1.

\textbf{Theorem 3.5.} Suppose $\A$ is approximately weakly
amenable and $\A$ is an ideal in $(\A^{**},\Box)$.
Then $\A$ is approximately $(2n-1)$-weakly amenable for every $n\in\Bbb N$.

{\it Proof.}  For $n\in Z^{+}$, we regard $\A^{(2n+2)}$ as the second dual
of $(\A^{(2n)},\Box)$ taken with the first Arens product $\Box$, and,
for $m\leq n$, we regard $(\A^{(2m)},\Box)$ as a subalgebra of $(\A^{(2n)},\Box)$.

Fix $n\in N.$ Since $\A$ is an ideal in $(\A^{**},\Box)$ then for each
$a\in\A$, the operators $L_{a}$ and $R_{a}$ on $\A$ are weakly
compact, and so the operators $L_{a}^{(2n)}$ and $R_{a}^{(2n)}$ are
weakly compact on $(\A^{(2n)},\Box)$. Thus
$a.\phi$,$\phi.a\in\A^{(2n-2)}$ for $a\in\A$ and $\phi\in\A^{(2n)}$.
Further, $a_{1}\ldots a_{n}.\phi$ and $\phi.a_{1}\ldots a_{n}$
belong to $\A$ for $a_{1},\ldots,a_{n} \in\A$ and $\phi\in\A^{(2n)}$.
Let the map $P:\A^{(2n-1)}\To\A^{*}$
be the natural projection. Then 
$$
\A^{(2n-1)}=\A^{*}\oplus \A^{\perp}
$$
as  $\A$-modules.

Suppose $D\in Z^{1}(\A,\A^{(2n-1)})$. Then there are derivation
$D_{1}\in Z^{1}(\A,\A^{*})$ and $D_{2}\in Z^{1}(\A,\A^{\perp})$ such that
$$Da=D_{1}a+D_{2}a\quad (a\in\A).$$
By assumption $D_{1}$ is approximately inner, and so it suffices for
the result to show that $D_{2}$ is approximately inner.

Let a,b$\in\A^{[n]}$where the later algebra is the n-th power of $\A$. For each $\phi\in\A^{(2n)}$ we have $\phi .a,\ b.\phi\in\A$ as we saw in the first paragraph and hence
$$
\langle D_{2}(ab),\phi\rangle=\langle D_{2}a.b,\phi\rangle+\langle a.D_{2}b,\phi\rangle=\langle D_{2}a,b.\phi\rangle+\langle D_{2}b,\phi.a\rangle=0
$$
So $D_{2}(ab)=0$. It follows that $D_{2}|_{\A^{[2n]}}=0$. By
Proposition 2.1 (i), $\A^{[2n]}$ is dense in $\A$, and so $D_{2}$=0. Therefore the
result follows. $\qed$

In the part (i) of the following Theorem we extend [2, Proposition 2.8.59(iii)] from $n=1$ to arbitrary $n\in\Bbb N$. In the next Theorem as it was pointed out in the remarks preceding Theorem 2.4, we consider $\A^{(2n)}$ as an $\A^{(2n)}$-module and then we equip $\A^{(2n+1)}$ with the usual dual module structure. 

\textbf{Theorem 3.6.} Let $n\in\Bbb N$ and  $B(\A,\A^{(2n-1)})=W(\A,\A^{(2n-1)})$. 

(i) Suppose that $(\A^{**},\Box)$ is $(2n-1)$-weakly amenable. Then $\A$ is
$(2n-1)$-weakly amenable.

(ii) Suppose $(\A^{**},\Box)$ is approximately $(2n-1)$-weakly amenable.
Then $\A$ is approximately $(2n-1)$-weakly amenable.

{\it Proof.} (i) Suppose $D\in Z^{1}(\A,\A^{(2n-1)})$. Take $\phi,\psi \in\A^{**}$, $\Gamma\in
\A^{(2n)}$ and bounded nets $(a_{\al})$ and $(b_{\beta})$ in $\A$ with
$$a_{\al}\To \phi\quad \text{and} \quad b_{\beta}\To \psi\quad \text{in}\quad \sigma(\A^{**},\A^*).
$$
Then
$$D^{**}(\phi\Box\psi)=\lim_{\al}\lim_{\beta}(D(a_{\al}).b_{\beta}+a_{\al}.D(b_{\beta})).$$
On the other hand we have

$$
\aligned \lim_{\al}\lim_{\beta}\langle D(a_{\al}).b_{\beta},\Gamma\rangle&=
\lim_{\al}\lim_{\beta}\langle D(a_{\al}),b_{\beta}\Box\Gamma\rangle\\
&=\lim_{\al}\langle D(a_{\al}),\psi\Box\Gamma\rangle=\langle D^{**}(\phi),\psi\Box\Gamma\rangle
=\langle  D^{**}(\phi).\psi,\Gamma\rangle .\endaligned
$$
 Our  assumption implies that $D\in W(\A,\A^{(2n-1)})$ and hence by [2, Theorem A.3.56] $D^{**}(\A^{**})\subseteq\A^{(2n-1)}$. Also the map 
$$
(A^{**},\sigma(A^{**}, A^{*}))\To (A^{(2n)},\sigma(A^{(2n)},A^{(2n-1)}),\quad \phi\mapsto\Gamma\Box \phi
$$ 
is continuous. So

$$
\aligned\lim_{\alpha}\lim_{\beta}\langle a_{\alpha}.D(b_{\beta}),\Gamma\rangle&=\lim_{\alpha}\lim_{\beta}\langle D(b_{\beta}),\Gamma.a_{\alpha}\rangle
=\lim_{\alpha}\langle D^{**}(\psi),\Gamma.a_{\alpha}\rangle\\
&=\lim_{\alpha}\langle D^{**}(\psi),\Gamma.a_{\alpha}\rangle=\langle D^{**}(\psi),\Gamma\Box\phi\rangle=\langle \phi.D^{**}(\psi),\Gamma\rangle.\endaligned
$$

Thus $D^{**}\in Z^{1}(\A^{**},\A^{(2n+1)})$ and so since $\A^{**}$ is
$(2n-1)$-weakly amenable, there exists $\La\in\A^{(2n+1)}$ such
that
$$D^{**}(\phi)=\phi.\La-\La.\phi\quad (\phi\in\A^{**}).$$
If  $P:\A^{(2n+1)}\To\A^{(2n-1)}$ is the natural
projection, then
$$
D(a)=a.P(\La)-P(\La).a\quad (a\in\A)
$$
Therefore $\A$ is $(2n-1)$-weakly amenable.

(ii) Suppose $D\in Z^{1}(\A,\A^{(2n-1)})$. As in part (i) $D^{**}\in
Z^{1}(\A^{**},\A^{(2n+1)})$. Thus there exists a net
$(\La_{\al})\subseteq\A^{2n+1}$ such that
$$
D^{**}(\phi)=\lim_{\al}(\phi.\La_{\al}-\La_{\al}.\phi)\quad (\phi\in\A^{**}).
$$
If $P:\A^{(2n+1)}\To\A^{(2n-1)}$ is the natural
projection, then
$$
D(a)=\lim_{\al}(a.P(\La_{\al})-P(\La_{\al}).a)\quad (a\in\A).
$$
Therefore $\A$ is approximately $(2n-1)$-weakly amenable. $\qed$

\textbf{Example 3.7.} An approximately $n$-weakly amenable Banach algebra which is not $n$-weakly amenable.

For each $n\in\Bbb N$ as in [9, Example 6.2] equip $M_{2^n}$ with the $\ell^2$ norm and let $\A_n$ be its unitization. If $\A=c_0(\A_n)$ then as it was shown in [9], $\A$ is approximately amenable but it is not weakly amenable. Therefore $\A$ is approximately $n$-weakly amenable for every $n\in\Bbb N$. However applying [3, Proposition 1.2] we conclude that it is not $(2n-1)$-weakly amenable for any $n\in\Bbb N$.    

\section{APPROXIMATE CYCLIC AMENABILITY}

Recall that a derivation $D:\A\To \A^{*}$ is called cyclic if $\langle Da,b\rangle+\langle a,Db\rangle=0$, $a,b\in \A$ and $\A$ is called cyclic amenable if every cyclic derivation $D:\A\To \A^{*}$ is inner. The natural approximate version of this concept is as follows. 

\textbf{Definition.} $\A$ is called approximately cyclic amenable if  every cyclic derivation $D:\A\To \A^{*}$ is approximately inner.

\textbf{Proposition 4.1.} Suppose that $\A$ is non-unital. Then $\A$ is approximately cyclic amenable if and only if $\A^{\#}$ is 
approximately cyclic amenable.

{\it Proof.}   Suppose $\A$ is approximately cyclic amenable and $D:\A^\#\To(\A^\#)^*$ is a cyclic derivation. Using the terminology and identities of Remark 3.2, we have $D(e)=0$ and one can observe that there exists a bounded   derivation $\wt D:\A\To\A^*$ and a bounded linear functional $\La\in\A^*$ such that $D(a)=\langle\La ,a\rangle e^*+\wt D(a)$, $a\in\A$. Since $D$ is cyclic, then for every $a,b\in\A$ we have
$$
0=\langle Da, b\rangle +\langle Db,a\rangle=\langle\wt Da, b\rangle +\langle\wt Db,a\rangle .
$$
So $\wt D$ is a cyclic derivation and hence by assumption there is a net $\{\La_\al\}$ in $\A^*$ such that $\wt D(a)=\lim_\al(a\La_\al-\La_\al a),\quad (a\in\A)$. Using the cyclic identity for both of $D$ and $\wt D$, for all $a,b\in\A$  we get 
$$
0=\langle Da,e+b\rangle +\langle D(e+b), a\rangle=\langle\La ,a\rangle + \langle\wt Da, b\rangle +\langle\wt Db,a\rangle .
$$
Thus $\La=0$ and hence 
$$
D(re+a)=\wt D(a)=\lim_\al(a\La_\al-\La_\al a)=\lim_\al[(re+a)\La_\al-\La_\al (re+a)].
$$
Therefore $\A^\#$ is approximately cyclic amenable. 

Conversely, suppose $\A^\#$ is approximately cyclic amenable and
$D:\A\To\A^{*}$ is a cyclic derivation. As in the proof of [12, Proposition 2.1]  
$$
\wt{D}:\A^{\#}\To(\A^{\#})^{*},\quad re+a\longmapsto Da,
$$ 
is a bounded derivation extending $D$; It is easy to see that $\wt D$ satisfies the cyclic identity. So there exists a net
$\{\La_\al\}$ in $\A^*$ and a net $r_\al$ in $\Bbb C$ such that
$$
\wt{D}(a+r)=\lim_\al[(a+r)(\La_\al+r_\al)-(\La_\al+r_\al)(a+r)]\quad (a\in\A,\quad r\in \Bbb C ).
$$
Therefore
$$
Da=\lim_\al (a.\La_\la-\La\al.a)\quad (a\in\A). \qed
$$
 
\textbf{Proposition 4.2.} Let $I$ be a closed ideal in $\A$.

(i) Suppose that $\A/I$ is approximately cyclic amenable. Then $I$ has the approximate trace extension property. 

(ii) Suppose that $\A$ is cyclic amenable and $I$ has the approximate trace extension property.
Then $\A/I$ is approximately cyclic amenable.

(iii) If $\A/I$ is approximately cyclic amenable, $\overline{I^2}=I$, and $I$ is cyclic amenable, then $\A$ is approximately cyclic amenable.

(iv) Suppose that $\A$ is approximately cyclic amenable and $I$ has a quasi-central bounded approximate identity for $\A$. Then $I$ is approximately cyclic amenable.

{\it Proof.} All parts can be proved with the same argument of Proposition 2.2, except that in the part (iii) the identity  $\overline{I^2}=I$, which is needed to complete the proof, does not follow from cyclic amenability of $I$ but it is a part of the assumptions. $\qed$

\textbf{Example 4.3.} An approximately cyclic amenable Banach algebra which is not cyclic amenable.

Let $\A$ be as in Example 3.7. Then $\A$ is approximately cyclic amenable, since it is approximately amenable. Let 
$$
P_1=\bmatrix 0&-1\\ 1&0\endbmatrix
$$  
and inductively define 
$$
P_{n+1}=\bmatrix 0&-P_n\\ P_n& 0\endbmatrix .
$$
Then define a derivation $D:\A\To\A^*=\ell^1(\A_n^*)$ by $D((x_n))=(\frac{1}{n^2}\delta_{P_n}(x_n))$. Then using the identity $\delta_P(B)=PB^t-B^tB$, one can observe that $D$ is cyclic. However as it was pointed out in [9, Example 6.2], $D$ is not inner.

\textbf{Example 4.4.} One might ask whether the condition $\overline{I^2}=I$ is necessary in Proposition 4.2(iii). The following example 
shows that this condition cannot be removed. 

Let $\A=\Bbb C^2$ with zero product and $I=\Bbb C\oplus 0$.  Then by [12, Example 2.5] $I$ and $\A/I$ are cyclic amenable but $\A$ is not cyclic amenable. Since for commutative Banach algebras, the two notions of cyclic amenability and approximate cyclic amenability coincide, then we see that $I$ and $\A/I$ are approximately cyclic amenable but $\A$ is not. Moreover $\overline{I^2}\neq I$.

\textbf{Remark 4.5.} In the Propositions 2.2 and 4.2 we observed that if we impose certain conditions to a closed ideal $I$ of $\A$, then approximate weak [respectively cyclic] amenability of $I$ and $\A/ I$ follows from that of $\A$. For subalgebras the situation is different, as the following proposition shows.

\textbf{Proposition 4.6.}  Let $B$ be a closed subalgebra of $\A$ and $I$ be a closed ideal of $\A$ such that $A=B\oplus I$. If $\A$ approximately weakly [respectively cyclic] amenable then so is $B$.

{\it Proof}. Suppose $\A$ is approximately weakly  amenable. Let $P:\A\To B$ be the natural projection  and $D:B\To B^{*}$ be a bounded derivation. By [16, Lemma 2.2] $P^{*}DP:\A\To A^{*}$ is a bounded derivation. So there exists a net $(\La_{\al})\subseteq A^{*}$ such that
$$
P^{*}DP(a)=\lim_{\al}(a.\La_{\al}-\La_{\al}.a)\quad (a\in A).
$$
For every $\al$ define $\la_{\al}=\La_{\al}|_{B}$. Observe that for every $b\in B$ and every $\al$ we have $P^{*}DP(b)|_B=Db$ and $(b.\La_\al -\La_\al .b)|_B=(b.\la_\al -\la_\al .b)$. So
$$
Db=P^{*}DP(b)|_B=\lim_{\al}(b.\La_{\al}-\La_{\al}.b))|_{B}=\lim_{\al}(b.\la_{\al}-\la_{\al}.b).
$$
Therefore $B$ is approximately weakly amenable. 

Now if $D:B\To B^{*}$ is a cyclic derivation, then observe that $P^*DP$ satisfies the cyclic identity. So by applying the same argument as above,  we obtain the statement for cyclic amenability. $\qed$ 

{\bf Acknowledgement.} The authors would like to express their sincere thanks to Professor F. Ghahramani for his valuable comments.

\vspace{20mm}

\noindent Department of Mathematics, Faculty of Sciences, Shiraz University, Shiraz 71454, Iran\\ \quad  esslamz@shirazu.ac.ir

\noindent Department of Mathematics, Karaj Branch, Islamic Azad University, Karaj, Iran\\  \quad
shoujaei@kiau.ac.ir

\end{document}